\documentclass[
leqno,          
11pt,           
a4paper         
]{amsart}

\usepackage[colorlinks,citecolor=blue,urlcolor=blue]{hyperref}          
\usepackage{courier}                                                    
\usepackage{amssymb, amsmath, amsfonts, amsthm}                         
\usepackage{mathtools, cite, enumerate, rotating, framed, lipsum}       
\usepackage{fancybox, bbm}                                              
\usepackage[dvipsnames]{xcolor}                                         
\usepackage[polish, english]{babel}                                     
\usepackage[T1]{fontenc}                                                
\usepackage[cp1250]{inputenc}                                           
\usepackage[normalem]{ulem}
\usepackage{yfonts}
\usepackage{mathrsfs}

\DeclareMathAlphabet{\mathpzc}{OT1}{pzc}{m}{it}
 \numberwithin{equation}{section}                        

\counterwithin{thmcount}{section}
\newcommand{\thmcount}{thmcount}                 
\newcounter{specialcounter}

\newtheorem{Thm}[\thmcount]{Theorem}
\newtheorem{Sthm}[specialcounter]{Theorem}
\newtheorem{Cor}[\thmcount]{Corollary}
\newtheorem{Lem}[\thmcount]{Lemma}
\newtheorem{Prop}[\thmcount]{Proposition}
\newtheorem{Rem}[\thmcount]{Remark}
\newtheorem{Defn}[\thmcount]{Definition}
\newtheorem{Ex}[\thmcount]{Example}
\newtheorem{Asu}[\thmcount]{Assumption}
\newtheorem{Sol}[\thmcount]{Solution}
\newtheorem*{Thmx}{Theorem}
\newtheorem*{Corx}{Corollary}
\newtheorem*{Lemx}{Lemma}
\newtheorem*{Propx}{Proposition}
\newtheorem*{Remx}{Remark}
\newtheorem*{Defnx}{Definition}
\newtheorem*{Exx}{Example}
\newtheorem*{Asux}{Assumption}
\newtheorem*{Solx}{Solution}

\newcommand \eq[1]{\begin{equation} #1 \end{equation}}
\newcommand \eqx[1]{\begin{equation*}  #1 \end{equation*}}
\newcommand \al[1]{\begin{align} #1 \end{align}}
\newcommand \alx[1]{\begin{align*}  #1 \end{align*}}
\renewcommand \sp[1]{\begin{equation} \begin{split} #1 \end{split} \end{equation}}
\newcommand \spx[1]{\begin{equation*} \begin{split} #1 \end{split} \end{equation*}}
\newcommand \en[1]{\begin{enumerate}  #1 \end{enumerate}}

\newcommand{\thm}[2]{\begin{Thm} \label{#1} #2 \end{Thm}}

\newcommand{\lem}[2]{\begin{Lem} \label{#1} #2 \end{Lem}}

\newcommand{\prop}[2]{\begin{Prop} \label{#1} #2 \end{Prop}}

\newcommand{\defn}[2]{\begin{Defn} \label{#1} #2 \end{Defn}}

\newcommand{\pr}[1]{\begin{proof} #1 \end{proof}}

\newcounter{comcount}

\renewcommand{\a}{\alpha}                
                 
        \newcommand{\la}{\lambda}

\newcommand{\EE}{\mathbb{E}}

\newcommand{\NN}{\mathbb{N}}  
  
\newcommand{\PP}{\mathbb{P}}  
  
\newcommand{\RR}{\mathbb{R}}

\newcommand{\supp}{\mathrm{supp}}
\newcommand{\8}{\infty}

\newcommand{\Rd}{{\RR^d}}

\newcommand{\wt}[1]{\widetilde{#1}}

\newcommand{\st}{\sqrt{t}}

\newcommand{\nat}{\mathbb{N}}

\newcommand{\abs}[1]{\left| #1 \right|}
\newcommand{\set}[1]{\left\{ #1 \right\}}
\newcommand{\norm}[1]{\left\| #1 \right\|}

\newcommand{\eee}[1]{\left( #1 \right)}

\newcommand{\wh}[1]{\widehat{#1}}

\newcommand{\sleq}{\lesssim}
\newcommand{\sgeq}{\gtrsim}

\begin{document}

\author[M. Dymowski \and M. Preisner \and A. Sikora]{Micha\l \ Dymowski \and Marcin Preisner \and Adam Sikora}

\address{Micha\l \ Dymowski, Instytut Matematyczny, Uniwersytet Wroc\l awski, \ pl. Grunwaldzki 2, 50-384 Wroc\l aw, Poland }
\email{michal.dymowski@uwr.edu.pl}

\address{Marcin Preisner, Instytut Matematyczny, Uniwersytet Wroc\l awski, \ pl. Grunwaldzki 2, 50-384 Wroc\l aw, Poland }
\email{marcin.preisner@uwr.edu.pl}

\address{
	Adam Sikora, Department of Mathematics, Macquarie University, NSW 2109, Australia}
\email[Corresponding author]{adam.sikora@mq.edu.au }

\title{Harmonic functions for Bessel operators}


\begin{abstract}
We verify the continuity of the Riesz transform from the operator related Hardy space to $L^1$ - 
Lebesgue space of integrable functions. For the standard Euclidean Laplace operator, this is a classical result that plays a significant role in harmonic analysis and theory of singular integral operators. Here, we consider a one-dimensional model of manifolds with ends and external Dirichlet boundary operators. This setting extends the work of Hassell and the third author. Specifically, we examine the real line with the measure $|x|^{n-1}dx$ leading to various versions of Bessel operators. For integer $n$, this mimics the measure on Euclidean $n$-dimensional space and the obtained results are expected to provide good predictions for a class of Riemannian manifolds with Euclidean ends. 
\end{abstract}

\maketitle
\section{Introduction}\label{sec1}
Suppose $d>0$ is a fixed parameter. Harmonic analysis related to the classical Bessel operator
\eq{\label{Bessel}
L_B f(x) = -f''(x) - \frac{d-1}{x}f'(x)
}
considered on the space $(0,\8)$ with the measure $x^{d-1}\, dx$ is a wide and deeply studied subject.
Since the publication of the seminal paper by Muckenhoupt and Stein \cite{Muckenhoupt_Stein}, this class of operators has garnered considerable attention.

If $d \in \NN$, then one can identify the Bessel operator $L_B$ with the usual Laplacian on $\Rd$ considered on radial functions. In other words, the Bessel process associated with $L_B$ is the length of the Brownian motion on $\Rd$. Harmonic analysis related to Bessel operators was studied by many authors, see e.g. \cite{BDT_d'Analyse, DP_Monats, Bessel12, Dziubanski_Bessel_Dunkl, Kania_Preisner, Bessel13, Bessel14}. 

In this note we shall consider $d\in (2,\infty)$ and the space $$\wt{X} = (-\8,-1] \cup [1,\8),$$
where the points $x=-1$ and $x=1$ are glued together, see Section \ref{sec3} below for precise description. We shall consider the Bessel process on $\wt{X}$. Every time the process arrives at $x=-1=1$ it chooses randomly the direction and can travel between the half-lines $(-\8,-1)$ and $(1,\8)$. For $d\in\nat_{+}$ this is a model of what happens with the process on manifolds of the form $\Rd \# \Rd$ (a connected sum of two copies of $\Rd$). Such a model was first considered by Hassell and the third author in \cite{Sikora_Hassell_1D}. 

The results obtained in  \cite{Sikora_Hassell_1D} are intriguingly related to several significant multidimensional studies on manifolds with ends or external Neumann or Dirichlet conditions (see, for example, \cite{KVZ, RiJ, Dangyang}). Our objective is to explore the anticipated relationship between the Riesz transform and Hardy spaces. Additionally, we present straightforward calculations for the corresponding heat kernel estimates. Studying these ideas in the described setting is more manageable than on manifolds with ends and provides a basis for proving expected results in this context.

We consider a metric-measure space $(\wt{X}, \rho, \mu)$, where $\rho$ is the Euclidean metric on both half-lines glued at the point $x=-1=1$ and $\mu$ is the measure with the density $|x|^{d-1} dx$. We shall denote by $\wt{L}$ the self-adjoint operator that corresponds to the Bessel process on $\wt{X}$, see Section~\ref{sec3} for details. Let $T_t = \exp(-t\wt{L})$ be the corresponding semigroup.

\defn{def_Hardy}{
For $f\in L^2(\wt{X})$ define
\begin{equation*}
    \label{def_Hardy_L}
 \norm{f}_{H^1_{\wt{L}}(\wt{X}) }:=\norm{\sup_{t>0}\abs{T_tf(x)} }_{L^1(\wt{X})}.
\end{equation*}
Following \cite{Hofmann_Memoirs}, the Hardy space $H^1_{\wt{L}}(\wt{X})$  is defined as the completion (in the norm $\norm{f}_{H^1_{\wt{L}}(\wt{X}) }$) of the set
$\set{f\in L^2(\wt{X}) \, : \, \norm{f}_{H^1_{\wt{L}}(\wt{X}) }< \8}
$.
}

In order to state this characterization for $d>2$ set
$$
h_+(x) = 
\begin{cases}
    1-\frac{1}{2}x^{2-d} & x\geq 1\\
    \frac{1}{2} |x|^{2-d} & x\leq -1
\end{cases}, \qquad h_-(x) = 
\begin{cases}
    \frac{1}{2}x^{2-d} & x\geq 1\\
    1-\frac{1}{2} |x|^{2-d} & x\leq -1
\end{cases}.$$
Notice that both functions $h_-,h_+$ are bounded and (locally weakly) $\wt{L}$-harmonic in the following sense: for $g\in C_c^\8(\wt{X})$ we have
 $\int_{\wt{X}} h_+'(x) g'(x)\, d\mu(x) =\int_{\wt{X}} h_-'(x) g'(x)\, d\mu(x) =0.$

In the corresponding Markov process, $h_+(x)$ represents the probability that the process, starting at $x$, will eventually escape to
$+\8$ as time progresses. Conversely, $h_- = 1-h_+$  denotes the probability that the process will escape via the opposite end, $-\8$.


\begin{Defn}
    \label{def_atomic}
We say that a function $a$ on $\wt{X}$ is an $(\mu,h_+,h_-)$-atom if there exists a ball $B$ in $\wt{X}$ such that:
\eq{\label{atoms_0}
\supp \, a \subseteq B, \quad \norm{a}_{L^2(X,\mu)} \leq \mu(B)^{-1/2}, \quad \int a h_+ \, d\mu = \int a h_- \, d\mu = 0.
}
The atomic Hardy space $H^1_{at}(\mu,h_+,h_-)$ is defined in the following way. A function $f$ belongs to $H^1_{at}(\mu,h_+,h_-)$ if $f(x) = \sum_k \la_k a_k(x)$, where $\sum_k |\la_k|<\8$ and $a_k$ are $(\mu,h_+, h_-)$-atoms. Then
\begin{equation}
    \label{eq_atomic_norm}
    \norm{f}_{H^1_{at}(\mu, h_+,h_-)} := \inf \sum_k |\la_k|,
\end{equation}
where the infimum is taken over all the representations as above.
\end{Defn}

By a standard argument the space $H^1_{at}(\mu, h_+,h_-)$ is a Banach subspace of $L^1(\wt{X}, \mu)$. The following proposition characterizes the space $H^1(L)$ in terms of atomic decomposition. 
In classical harmonic analysis, specifically for the Laplacian $\Delta$ on $\Rd$, such theorem was established in \cite{Coifman_Studia} and \cite{Latter_Studia}. For further generalizations, see \cite{Hofmann_Memoirs, DP_Annali, Preisner_Sikora_Yan, Preisner_Sikora}.

\prop{thmA}{
Let $d>2$ and $\wt{L}$ be the Bessel operator on $\wt{X}$. Then the spaces $H^1_{\wt{L}}(\wt{X})$ and $H^1_{at}(\mu,h_+,h_-)$ coincide and there exists $C>0$ such that
$$
C^{-1} \norm{f}_{H^1_{\wt{L}}(\wt{X})} \leq \norm{f}_{H^1_{at}(\mu,h_+,h_-)}\leq C \norm{f}_{H^1_{\wt{L}}(\wt{X})}.
$$
}
Proposition \ref{thmA} is a consequence of \cite{Preisner_Sikora_Yan} and \cite{Preisner_Sikora}, see Section \ref{ssec3.3} below. Recall that the Riesz transform $\wt{R}$ related to $\wt{L}$ is given by $\wt{R} = \partial_x \wt{L}^{-1/2}$. It was proved in \cite{Sikora_Hassell_1D} that $\wt{R}$ is bounded on $L^p(\wt{X}, \mu)$ for $p\in(1,d)$ and unbounded for $p\geq d$. One of the main results in this paper is the following end-point at $p=1$ estimate for the Riesz transform $\wt{R}$.
\thm{thmB}{
Suppose that $d>2$. The Riesz transform $\wt{R}$ related to the Bessel operator on $\wt{L}$ is bounded from $H^1_{\wt{L}}(\wt{X})$ to $L^1(\wt{X},\mu)$, i.e. there exists $C>0$ such that
$$
\norm{\wt{R} f}_{L^1(\wt{X},\mu)}\leq C \norm{f}_{H^1_{\wt{L}}(\wt{X})}.$$
}

As we will demonstrate, a key and intriguing aspect is that atoms in $H^1_{at}(\mu, h_+,h_-)$ necessarily satisfy the cancellation conditions with respect to both functions $h_+$ and $h_-$. The proof of Theorem \ref{thmB} is provided in Section \ref{ssec3.4} below. It is important to note that the operator $\wt{R}$ is not bounded on the classical Coifman-Weiss atomic Hardy space, as discussed in Remark~\ref{remark} below.

The following statement describes behaviour of the heat kernel $T_t(x,y)$. It is closely related to the
estimates obtained by Grigor'yan and Saloff-Coste  
\cite{Grigoryan_Saloff-Coste}. {However, the proof described below is noteworthy for its conciseness.} 
\thm{thmC}{
Assume that $d>2$ and let $T_t (x,y)$ be the kernel of the semigroup $T_t = \exp(-t\wt{L})$. Then there exists $C,c_1,c_2>0$ such that for $x,y\in \wt{X}$ we have:
\al{ \label{ULmixedBounds0} \frac{C^{-1}}{\mu(B(x,\sqrt{t}))} e^{-\frac{\rho(x,y)^2}{c_1 t}} \leq T_t(x,y) &\leq  \frac{C}{\mu(B(x,\sqrt{t}))} e^{-\frac{\rho(x,y)^2}{c_2t}}, &x\cdot y >0\\
\label{ULmixedBounds}
C^{-1} \frac{|x|^{2-d} + |y|^{2-d}}{\mu(B(x,\sqrt{t}))} e^{-\frac{\rho(x,y)^2}{c_1 t}} \leq T_t(x,y) &\leq C \frac{|x|^{2-d} + |y|^{2-d}}{\mu(B(x,\sqrt{t}))} e^{-\frac{\rho(x,y)^2}{c_2t}},  &x\cdot y <0.
}
}

{\bf Organisation of the paper.} In Section \ref{sec2} we study the Bessel operators on $[1,\8)$ with either Neumann or Dirichlet boundary conditions. In particular, we prove Theorem \ref{thm_riesz_halfline} which is the main tool in the proof of Theorem \ref{thmB}. The proofs of  Proposition \ref{thmA}, Theorem \ref{thmB}, and Theorem \ref{thmC} are given in Section \ref{sec3} below.

\section{Bessel operators on $[1,\8)$.}\label{sec2}

In this section we study two Bessel operators related to the differential operator \eqref{Bessel} on the space $X_+ := [1,\8)$ with the measure $d\mu(x) = x^{d-1}\, dx$. The first operator $L_N$ will be related to the Neumann boundary condition at the boundary $x=1$, whereas the second $L_D$ will have the Dirichlet boundary condition at $x=1$.

\begin{Defn}
\label{def_1}
Let $\Sigma_N=
\set{ f\in C^1_c(X_+) \,  : \, f'(1)=0}\subseteq L^2(X_+,d\mu).$ Define
$$Q(f,g)=
\int_{X_+} f'(x)g'(x)d\mu(x), \qquad f,g\in \Sigma_N .$$
Let $Q_N$ the closure of the quadratic form $Q$. The corresponding self-adjoint operator $L_N$ satisfies $Q_N(f,g) = \langle L_N f, g \rangle$. This operator is called the Bessel operator on $X_+$ with the Neumann boundary condition at $x=1$. Analogously, we define the  Bessel operator  $L_D$ with the Dirichlet boundary condition and its quadratic form using initially the core $\Sigma_D = \set{f \in C^1_c(X_+) : f(1)=0 }.$ 
\end{Defn}

If the parameter $d>0$ is an integer, then $L_N$ can be considered as the distance from zero of the standard Brownian motion that lives on the complement of the unit ball and is reflected at the boundary. The operator $L_D$ can be interpreted similarly, but with the process that is killed at the boundary.

Here we are interested mainly in the case $d>2$. Denote by $T_{t,N}$ and $T_{t,D}$ the semigroups generated by $L_N$ and $L_D$, respectively. The Hardy spaces $H^1_{L_N}(X_+)$ and  $H^1_{L_D}(X_+)$ are defined analogously as in Definition \ref{def_Hardy} with the semigroups $T_t^N$ and $T_t^D$, respectively. The Riesz transforms $R_N = \partial_x L_N^{-1/2}$ and $R_D = \partial_x L_D^{-1/2}$ can be expressed by the resolvents:
\begin{equation}
\label{Riesz_resolvent}
R_N= \partial_x \int_0^{\infty}\frac{2}{\pi}  \left(L_N + k^2\right)^{-1} \, dk, \quad R_D= \partial_x \int_0^{\infty}\frac{2}{\pi}  \left(L_D + k^2\right)^{-1} \, dk.
\end{equation}

Our main goal in this section is to prove the following boundedness result.
\begin{Thm}
 \label{thm_riesz_halfline}
Let $d>2$. Then:
\en{
\item $R_N$ is bounded from $H^1_{L_N}(X_+)$ to $L^1(X_+,\mu)$,
\item $R_D$ is bounded from $H^1_{L_D}(X_+)$ to $L^1(X_+,\mu)$.
}
\end{Thm}

The proof of Theorem \ref{thm_riesz_halfline} is given in Section \ref{ssec24}. It is based on estimates from Proposition~\ref{prop_kernel} below. These estimates are obtained by quite detailed analysis of the kernels $R_N(x,y)$ and $R_D(x,y)$ that are expressed by the related resolvent kernels, see \eqref{Riesz_resolvent}.

\subsection{Bessel operator on $[1,\8)$ with Neumann boundary condition.}\label{ssec21}
Let $d>0$. It is well known that there exists an integral kernel $T_{t,N}(x,y)$ associated with $T_{t,N}$ and positive constants $C,c_1,c_2$ such that
\eq{\label{gauss1}
C^{-1} \mu(B(x,\st))^{-1} \exp\eee{-\frac{|x-y|^2}{c_1 t}}\leq T_{t,N}(x,y) \leq C \mu(B(x,\st))^{-1} \exp\eee{-\frac{|x-y|^2}{c_2 t}},
}
for $x,y\in X_+$ and $t>0$. The unique (up to a multiplicative constant) bounded harmonic function for $L_N$ is $h_N \equiv 1$.

Denote by $H^1_{at}(X_+)$ the Coifman-Weiss atomic Hardy space, see \cite{CoifmanWeiss_BullAMS}. Recall that a function $a$ on $X_+$ is a Coifman-Weiss atom if there exists a ball $B$ in $X_+$ such that:
\eq{\label{Coifman-Weiss}
\supp \, a \subseteq B, \quad \norm{a}_{L^2(X_+,\mu)} \leq \mu(B)^{-1/2}, \quad \int_B a\, d\mu = 0.
}
The Coifman-Weiss Hardy space $H^1_{at}(X_+)$ is defined analogously as in Definition \ref{def_atomic}, using atoms as above. The following theorem states that the space $H^1_{L_N}(X_+)$ can be characterized in terms of atoms.
\begin{Thm}
    \label{thm_Hardy_Neu}
    Let $d>0$ and $L_N$ be the Bessel operator on $[1,\8)$ with Neumann boundary condition. Then the spaces $H^1_{L_N}(X_+)$ and $H^1_{at}(X_+)$ coincide and have comparable norms.
\end{Thm}
Theorem \ref{thm_Hardy_Neu} follows from  \cite[Theorem 2]{DP_Annali}, see also \cite{Preisner_Sikora_Yan} and \cite{Preisner_Sikora}.

The Riesz transform $R_N = \partial_x L_N^{-1/2}$ was studied in \cite{Sikora_Hassell_1D}, where it was proved that $R_N$ is bounded on $L^p(X_+,\mu)$ for $p\in(1,\8)$. Theorem~\ref{thm_riesz_halfline} is an end-point result for $R_N$.

\subsection{Bessel operator on $[1,\8)$ with Dirichlet boundary condition.}\label{ssec22}

The analysis of $L_D$ is a bit different. The kernel $T_{t,D}(x,y)$ satisfies the upper estimate from \eqref{gauss1}, but not the  lower one. From now on we consider $d>2$. The harmonic function for $L_D$ is
\begin{equation} \label{eq_hD}
    h_D(x) = 1 - x^{2-d}.
\end{equation}
Denote by $\mu_{h_D^2}$ the measure on $X_+$ with the density $h^2_D(x) \, d\mu(x)$. The following estimates follow from the results of Gyrya and Saloff-Coste \cite{Gyrya_Saloff-Coste}. There exist $c_1,c_2>0$ such that for $t>0$ and $x,y\in X_+$ we have
\eq{\label{gauss2}
 \mu_{h_D^2}(B(x,\st))^{-1} \exp\eee{-\frac{|x-y|^2}{c_1 t}}\sleq \frac{T_{t,D}(x,y)}{h_D(x)h_D(y)} \sleq \mu_{h_D^2}(B(x,\st))^{-1} \exp\eee{-\frac{|x-y|^2}{c_2 t}}.
}

To describe the Hardy space $H^1_{L_D}(X_+)$, following \cite{Preisner_Sikora}, define the $(X_+, h_D)$-atoms by the conditions:
\eq{\label{atom_Dir}
\supp\, a \subseteq B, \quad \norm{a}_{L^2(X_+,\mu)} \leq \mu(B)^{-1/2}, \quad \int_B a \,h_D \, d\mu =0,
}
compare \eqref{Coifman-Weiss}. The Hardy space $H^1_{at}(X_+, h_D)$ is defined as in Definition \ref{def_atomic}, with $(X_+, h_D)$-atoms. Similarly as in the case of  $H^1_{L_N}(X_+)$, the space $H^1_{L_D}(X_+)$ can be characterized in terms of atoms.

\begin{Thm}
    \label{thm_Hardy_Dir}
    Let $d>2$ and $L_D$ be the Bessel operator on $[1,\8)$ with the Dirichlet boundary condition at $x=1$. Then the spaces $H^1_{L_D}(X_+)$ and $H^1_{at}(X_+, h_D)$ coincide and have comparable norms.
\end{Thm}
Theorem \ref{thm_Hardy_Dir} is a consequence of \eqref{gauss2} and \cite[Theorem A]{Preisner_Sikora}.

The Riesz transform $R_D=\partial_x L_D^{-1/2}$ was studied in \cite{Sikora_Hassell_1D}. It was demonstrated that that $R_D$ is bounded on $L^p(X_+,\mu)$ for $p\in(1,d)$ and unbounded for $p\geq d$. An endpoint result for $p=d$ was subsequently established in \cite{Dangyang}.

\subsection{Preliminary estimates.}\label{ssec23}

In this section we shall study the singular integral kernels $R_N(x,y)$ and $R_D(x,y)$ associated with $R_N$ and $R_D$, respectively. Using the resolvent representation of the Riesz transform, the authors in \cite{Sikora_Hassell_1D} obtained the explicit formulas for $R_N(x,y)$ and $R_D(x,y)$ given in terms of the Bessel functions $I_\a$ and $K_\a$, see Proposition \ref{prop_Riesz_kernel}. We shall need the following estimates for the derivatives of these kernels.

\prop{prop_kernel}{
Let $d>2$. There exists $C>0$ such that for $x,y\in X_+$ we have:
 \al{\label{Riesz_est_N}
 \abs{ \partial_y R_{N} (x,y)} \leq  & Cx^{1-d}(x-y)^{-2},\\
  \label{Riesz_est_D}
 \abs{\partial_y \eee{\frac{R_D(x,\cdot)}{h_D(\cdot)}} (y)} \leq & C \frac{1}{h_D(y)} \frac{1}{|x-y|^2}x^{-1}\min(x^{2-d},y^{2-d}).
 }
}
The proof of Proposition \ref{prop_kernel} requires quite detailed analysis and is given below on page \pageref{pr_lem4}.  
For $z>0$ denote: 
\spx{
l(z) = z^{1-d/2} I_{d/2-1}(z),  \qquad k(z) = z^{1-d/2} K_{d/2-1}(z),
}
where $I_\a$ and $K_\a$ are the modified Bessel functions, see \cite{Abramowitz_Stegun}.  Following \cite{Sikora_Hassell_1D} denote $A = l/k$, $B=l'/k'$. Next, for $\la >0$ and $z>1$ set: 
\sp{\label{tau}
\tau_\la(z) = l(\la z) - A(\la) k (\la z), \qquad 
\psi_\la(z) &= l(\la z) - B(\la) k (\la z).
}

\prop{prop_Riesz_kernel}{\cite[Section 5]{Sikora_Hassell_1D}
The operators $R_{N}$ and $R_{D}$ are associated with the following kernels:
\eqx{
R_{N}(x,y) =
\begin{cases}
\int_0^\8  \la^{d-2} \psi_\la'(x) k(\la y) \, d\la & \text{if } \ 1<x\leq y\\
\int_0^\8 \la^{d-1} \psi_\la( y) k'(\la x)  \, d\la & \text{if } \ 1<y\leq x,
\end{cases}
}
\eqx{
R_{D}(x,y) =
\begin{cases}
\int_0^\8  \la^{d-2} \tau_\la'(x) k(\la y) \, d\la & \text{if } \ 1<x\leq y\\
\int_0^\8 \la^{d-1}  \tau_\la( y) k'(\la x)  \, d\la & \text{if } \ 1<y\leq x.
\end{cases}
}
}
The table below summarizes the asymptotics that are required for our discussion below. Since the asymptotics are slightly different in the cases $d\in (2,3)$ and $d\geq 3$ it is convenient to denote
\eq{\label{kappa}
\nu =
\begin{cases}
d-2 & \text{ if } d\in (2,3)\\
1 & \text{ if } d\geq 3
\end{cases}.
}

The following asymptotics are well known. We refer the reader to e.g. \cite{Abramowitz_Stegun}. Notice that we need here a slightly stronger estimates that the ones used in \cite{Sikora_Hassell_1D}. Denote $\wt{k}(z) = z^{d-2} k(z)$.

\setlength{\tabcolsep}{8pt}
\setlength{\arrayrulewidth}{0.3mm}
\renewcommand{\arraystretch}{1.4}

\begin{center}
\begin{tabular}{|c|c|c|c|}
  \hline
  function & sign  & $z\sim 0$ & $z \sim \8 $ \\
  \hline
 $l(z)$ & $ + $    & $c_l + O(z^2)$ & $\wt{c}_lz^{(1-d)/2} e^z\eee{1  + O(z^{-1})}$\\
 \hline
 $ l'(z) $ & $+$ & $ c_{l'} z \eee{1 + O(z^2)} $ & $ \wt{c}_lz^{(1-d)/2} e^z\eee{ 1  + O(z^{-1})} $ \\
 \hline
 $ k(z) $ & $+$ & $  c_k z^{2-d}  + d_k z^{2-d+\nu} + O(z^{3-d}) $ & $ \wt{c}_k z^{(1-d)/2} e^{-z}\eee{ 1  + O(z^{-1})} $ \\
 \hline
 $ k'(z) $ & $-$ & $ -c_{k'} z^{1-d}\eee{1 + O(z^{2})} $ & $ - \wt{c}_k z^{(1-d)/2} e^{-z}\eee{ 1  + O(z^{-1})} $ \\
 \hline
 $ A(z) $ & $+$ & $ \frac{c_l}{c_k} z^{d-2}\eee{1 + O(z^{\nu})} $ & $ \frac{\wt{c}_l}{\wt{c}_k}e^{2z} \eee{ 1  +
 O(z^{-1})} $ \\
 \hline
 $ B(z) $ & $-$ & $ -\frac{c_{l'}}{c_{k'}} z^{d}
 \eee{1+ O(z^2)} $ & $ -\frac{\wt{c}_l}{\wt{c}_k}e^{2z} \eee{ 1  + O(z^{-1})} $ \\
 \hline
 $ \wt{k}(z) $ & $+$ & $ c_k + d_k z^\nu + O(z) $ & $ \wt{c}_k z^{(d-3)/2} e^{-z}\eee{ 1  + O(z^{-1})} $ \\
 \hline
 $ \wt{k}'(z) $ &  & $ O(z^{\nu-1}) $ & $ -\wt{c}_kz^{(d-3)/2} e^{-z}\eee{ 1  + O(z^{-1})} $ \\
 \hline
\end{tabular}
\end{center}
\label{table_page}
The constants appearing in the above table satisfy: 
$$c_{k'}/c_k = d-2. \quad \wt{c}_{l} = (2\pi)^{-1/2}, \quad c_l = \Gamma(d/2)^{-1}, \quad \wt{c}_k = \sqrt{\pi/2}.$$
Recall that the harmonic function for $L_{D}$ is $h_D(x) = 1-x^{2-d}$.

\lem{tau_lemma}{
The function $\tau_\la(z)$, see \eqref{tau}, for $\la>0$ and $z>1$ has the following properties.
\en{
\item If $z>0$, then $\tau_\la (z), \tau'_\la(z) >0$.
\item If $\la z <1$, then
\al{
\label{tau_1}
\tau_\la (z) &=c_l h_D(z)\eee{1  -  \frac{d_k\la^\nu}{\wt{k}(\la)}} + O(\la (z-1)),\\
\label{tauprim_1}
\tau_\la'(z) &= c_lh_D'(z)\eee{1  -  \frac{d_k\la^\nu}{\wt{k}(\la)}} + O(\la),
}
where $d_k$ is the constant from the table above.
\item If $\la z >1$, then
\al{
\label{tau_2}
\tau_\la(z) &\leq C (\la z)^{(1-d)/2} \min(1, \la(z-1)) e^{\la z},\\
\label{tauprim_2}
\tau_\la'(z) &\leq C \la (\la z)^{(1-d)/2} e^{\la z}.
}
}
}

\pr{
{\bf 1.} Notice that $\tau_{\la}(z) =  k(\la z)^{-1} \eee{\frac{l}{k}(\la z) - \frac{l}{k}(\la)} >0 $
since $(l/k)' = k^{-2} (l'k-k'l) >0$, see \cite[Eq. (3.1)]{Sikora_Hassell_1D}, and $\la z > \la$. Then,
$$
\tau'_\la(x) = \la l'(\la z) - \la A(\la) k'(\la z)>0,
$$
since $l',A >0$ and $k'<0$.

{\bf 2.} Here we deal with $\la z<1$, so also $\la<1$. Write
\spx{
\tau_\la(z) = &c_l h_D(z)+ \underbrace{h_D(z)\eee{l(\la)-c_l}}_{A_1}+\underbrace{l(\la z)-l(\la)}_{A_2}\\
&+ \underbrace{(l(\la)-c_l) \eee{z^{2-d} - \frac{k(\la z)}{k(\la)} }}_{A_3}+\underbrace{c_l\eee{z^{2-d} - \frac{k(\la z)}{k(\la)} }}_{A_4}.
}
Recall that $h_D(z) \simeq (z-1)/z \leq (z-1)$. Therefore:
\alx{
\abs{A_1} &\sleq  h_D(z) \la^2 \sleq \la (z-1),\\
\abs{A_2} &= C\la (z-1) l'(\xi) \sleq \la(z-1).
}
where $\la<\xi<\la z <1$. Let $\rho$ be defined by $\wt{k}(t) = c_k + d_k t^\nu + \rho(t)$. Then $|\rho'(t)|\sleq t^{\nu-1} \leq C$. Next,
\sp{\label{eq43222}
 \frac{k(\la z)}{k(\la)} - z^{2-d} &= \frac{z^{2-d}}{\wt{k}(\la)}\eee{\wt{k}(\la z) - \wt{k}(\la)} \\
&= \frac{d_k z^{2-d}}{\wt{k}(\la) }\eee{\la^\nu(z^\nu -1) + \rho(\la z) - \rho(\la)}\\
&= \frac{d_k \la^\nu }{\wt{k}(\la) } z^{2-d}(z^{\nu} -1) + O(\la(z-1))
}
Consider first the case $d\in(2,3)$. Notice that $z^{2-d} (z^{\nu}-1) = h_D(z)$. Therefore,
\alx{
|A_3| &\sleq \la^{1+\nu} h_D(z) \leq \la(z-1),\\
A_4 &= -c_l h_D(z) \frac{d_k \la^\nu}{\wt{k}(\la)} + O(\la(z-1)).
}
If $d\geq 3$ the estimates are easier, since $\nu =1$ and $\la^\nu h_D(z) = O(\la(z-1))$. Then $|A_3+A_4| \sleq \la(z-1)$.

Now let us deal with $\tau_\la'(z) = \la l'(\la z) - \la A(\la) k'(\la z)$. Note that $|\la l'(\la z)| \sleq  \la^2 z \sleq \la$. Moreover,
\spx{
-\la A(\la) k'(\la z) =& 
\underbrace{\frac{c_l c_{k'}}{c_k} z^{1-d}\eee{1-
\frac{d_k\la^{\nu}}{\wt{k}(\la)}}}_{A_5} 
+ \underbrace{\frac{c_{k'}}{c_k}z^{1-d} \eee{c_k
\la^{2-d}A(\la)-c_l+\frac{d_k\la^{\nu}}{\wt{k}(\la)}}}_{A_6} \\
&-
\underbrace{\la A(\la)\eee{ c_{k'} (\la z)^{1-d}+k'(\la z)}}_{A_7}.
}
Recall that $c_{k'}/c_k = d-2$, so that $A_5 = c_l(d-2) z^{1-d}= c_l h'(z)$. Also, estimates from the table on page \pageref{table_page},
\alx{
|A_6| \sleq 
\abs{\eee{c_k\la^{2-d}A(\la)-c_l+\frac{c_l d_k\la^{\nu}}{\wt{k}(\la)}}}
= C\wt{k}^{-1}(\la)\abs{\eee{c_kl(\la)-c_l\wt{k}(\la) + c_l d_k\la^{\nu}}}
\sleq \la}
and
\alx{
|A_7| 
\sleq  \la \la^{d-2} (\la z)^{3-d} 
\sleq  \la (\la z) z^{2-d} 
\sleq  \la.
}
Thus, \eqref{tauprim_1} is proved.

{\bf 3.}  We deal with $\la z>1$. Notice that $A(\la) \leq C e^{2\la}$ holds for all $\la>0$ and $e^{\la(2-z)}\leq e^{\la z}$. Therefore, \eqref{tauprim_2} follows from the asymptotics of $l',k',A$, etc. Also, in the case when $\la(z-1)\geq 1/2$, the estimate \eqref{tau_2} can be verified in a  similar way.

In the remaining case, $\la (z -1)<1/2$, we have $\la = \la z - \la (z-1) > 1-1/2 = 1/2$. Denote $L(z) = l(z) z^{(d-1)/2} e^{-z}$ and $K(z) = k(z) z^{(d-1)/2} e^{z}$, so that $L(z) \simeq C$ and $K(z)\simeq C$  on $(1/2, \8)$. Now,
\eq{\label{eq456}\tau_\la(z) = l(\la) z^{(1-d)/2} e^{\la(z-1)} \underbrace{\eee{\frac{L(\la z)}{L(\la)} -1 +  e^{2\la(1-z)} \eee{1- \frac{K(\la z)}{K(\la)}}  + 1-e^{2\la(1-z)}  }}_{A_8}.
}
Using the mean value theorem and the fact that $L'(z)$ and $K'(z)$ are bounded on $(1/2,\8)$ we get that $|A_8| \sleq \la(z-1)$ and, consequently,
\sp{
\label{eq654}|\tau_\la(z)|&\sleq  \la^{(1-d)/2} e^{\la} z^{(1-d)/2} e^{\la(z-1)} \la(z-1) \sleq (\la z)^{(1-d)/2} e^{\la z} \la (z-1),}
which ends the proof of \eqref{tau_2}.
}

\lem{lem5}{
The function $\psi_\la(z) = l(\la z) - B(\la) k(\la z)$, $\la>0, z>1$, has the following properties:
\en{
\item $\psi_\la (z)>0$.
\item For $\la z <1$ we have:
\alx{
\psi_\la (z) &\leq C,\\
|\psi_\la'(z)| &\leq C\la^2 z .
}
\item For $\la z >1$ we have:
\alx{
\psi_\la(z) &\leq C (\la z)^{(1-d)/2} e^{\la z},\\
|\psi_\la'(x)| &\leq C \la (\la x)^{(1-d)/2} e^{\la x}.
}
}
}

\pr{[Sketch of the proof.]
The proof follows directly from the asymptotics from the table on page \pageref{table_page}. This proof is simpler than the proof of Lemma \ref{tau_lemma}, since we do not have to deal with subtle cancellations. The details are left to the reader.
}

\begin{proof}
[Proof of Proposition \ref{prop_kernel}]\label{pr_lem4}
We shall only prove \eqref{Riesz_est_D}. The proof of \eqref{Riesz_est_N} is similar and much simpler. In the estimates below we shall use Lemma \ref{tau_lemma} and asymptotic estimates for $l,k,l',k',$ etc. without further mention.

To prove \eqref{Riesz_est_D} set
\eq{\label{Kdef}
K(x,y) =  \abs{\partial_z \left.\eee{\frac{R_D(x,z)}{h_D(z)}} \right|_{z=y}}.
}

{\bf Case 1: $x\leq y$.} By Proposition \ref{prop_Riesz_kernel}
\alx{\label{eq2}
 K(x,y)&\leq  
 \int_0^\8 \la^{d-2} \tau_\la' (x) \abs{\partial_z\left.\eee{\frac{k(\la 
z)}{h_D(z)}}\right|_{z=y}}   \, d\la\\
 &\leq \int_0^\8 \la^{d-2} \tau_\la' (x) h_D(y)^{-2} V_\la(y)   \, d\la,
}
where
$
V_\la (y) = 
 |\la h_D(y) k'(\la y) - h_D'(u)  k(\la y)|.
$
Write
\spx{
\label{eq2}
K(x,y)\leq h^{-2}_D(y) \int_0^\8 \la^{d-2} \tau_\la'( x) V_\la(y)  \, d\la &= h^{-2}_D(y) \eee{\int_0^{1/y} ... + \int_{1/y}^{1/x} ... + \int_{1/x}^\8 ... }\\
&=I_1+I_2+I_3.
}
By Lemma \ref{tau_lemma},
\alx{
\tau'_\la(x)& \sleq  x^{1-d}(1+ \la^{\nu}) + \la 
\sleq x^{-1}, && \la x \leq 1,\\
\tau'_\la(x) &\sleq  \la (\la x)^{(1-d)/2} e^{\la x}, && \la x \geq 1,\\
V_\la(y)&\sleq  \la (\la y)^{1-d}+  y^{1-d} (\la y)^{2-d} \leq C \la  (\la y)^{1-d}, && \la y \leq 1,\\
V_\la(y) &\sleq  (\la y)^{(1-d)/2} e^{-\la u} (\la h_D(y) + h_D'(y)), && \la y \geq 1.}

Using the estimates above we arrive at:
\alx{
|I_1| 
&\sleq h_D^{-2}(y)  \int_0^{1/y} \la^{d-2} x^{-1} \la (\la y)^{1-d} \, d\la 
=  h_D^{-2}(y)x^{-1}y^{-d}, \\
|I_2|
&\sleq h_D^{-2}(y) \int_{1/y}^{1/x} \la^{d-2} x^{-1}   \la (\la y)^{(1-d)/2} e^{-\la y} \, d\la\\
&\sleq h_D^{-2}(y)  x^{-1} y^{(1-d)/2} \int_{1/y}^\8 \la^{(d-1)/2}e^{-\la y} \, d\la\\
&\sleq h_D^{-2}(y)   x^{-1} y^{-d}.
}
Recall that $h_D(y) \simeq (y-1)/y$ and $y> y-1 > y-x$. Hence,
$$|I_1|+|I_2| \sleq h(y)^{-1} x^{-1} \frac{y^{2-d}}{y(y-1)}  \sleq h(y)^{-1} x^{-1}\frac{y^{2-d}}{(y-x)^2}.$$
Notice that we have $h_D'(y)/h_D(y) \simeq y^{2-d}/(y-1) \leq 1/(y-x)$. Therefore,
\spx{
|I_3|&\sleq h_D^{-2}(y)  \int_{1/x}^\8  \la^{d-1} (\la x)^{(1-d)/2} e^{\la x} (\la y)^{(1-d)/2} e^{\la y} \eee{\la h_D(y) + h_D'(y)} \, d\la\\
&\sleq    h_D^{-1}(y)(xy)^{(1-d)/2} \int_{1/x}^\8 e^{-\la (y-x)} \frac{\la h_D(y) + h_D'(y)}{h_D(y)}\, d\la\\
&\sleq    h_D^{-1}(y)(xy)^{(1-d)/2} (y-x)^{-2} e^{-c y/x}\\
&\sleq   h_D^{-1}(y)(xy)^{(1-d)/2} (y-x)^{-2} (x/y)^{(d-3)/2} \sleq h_D^{-1}(y)x^{-1}\frac{y^{2-d}}{(y-x)^2}.
}
where we have used that $\int_{(y-x)/x}^\8  \la^M e^{-\la}\, d\la \sleq e^{-cy/x}$ for $M=0,1$. Thus  \eqref{Riesz_est_D} is proved. 

{\bf Case 2: $y\leq x$.}
By Proposition \ref{prop_Riesz_kernel},
\alx{\label{eq2}
 K(x,y)&\leq  
 \int_0^\8 \la^{d-1} |k' (\la x)| \abs{\partial_z \left.\eee{\frac{\tau_\la(z)}{h_D(z)}}\right|_{z=y}}   \, d\la\\
 &\leq \int_0^{1/x}... + \int_{1/x}^{1/y}...+ \int_{1/y}^{\8}... = J_1+J_2+J_3.
}

For $J_1$ and $J_2$ we have $\la y \leq 1$. Next, by \eqref{tau_1},
\spx{
\abs{\partial_y\frac{\tau_\la( y)}{h_D(y)}} 
=& h_D^{-2}(y) \abs{\tau'_{\la}(y)h_D(y) - \tau_{\la}(y)h'_D(y) } \\
=& h_D^{-2}(y) \left|\eee{c_lh_D'(z)\eee{1  -  \frac{d_k\la^\nu}{\wt{k}(\la)}} + O(\la)}h_D(y) \right.\\
&\left.- \eee{c_lh_D(z)\eee{1 
-  \frac{d_k\la^\nu}{\wt{k}(\la)}} + O(\la (y-1))}h'_D(y)\right| \\
\leq  &h_D^{-1}(y)\eee{ O(\la) + \frac{h'_D(y)}{h_D(y)}O(\la (y-1))}
\leq h_D^{-1}(y)O(\la).
}
Thus
\alx{
J_1 &\sleq   h_D^{-1}(y)\int_0^{1/x}\la^{d-1} (\la x)^{1-d} \la \, d\la 
\sleq h_D^{-1}(y)x^{-d-1},\\
J_2 
&\sleq h_D^{-1}(y) \int_{1/x}^{1/y}\la^{d-1} (\la x)^{(1-d)/2}e^{-\la x} \la \, d\la\\
&\sleq h_D^{-1}(y) x^{(1-d)/2} \int_0^\8 \la^{(d+1)/2}e^{-\la x}\, d\la 
= h_D^{-1}(y) x^{-d-1}.
 }

Recall that $h_d(y) \simeq (y-1) h_D'(y)$. For $\la y >1$, by Lemma \ref{tau_lemma},
\spx{\label{eq_West}
\abs{\partial_y\eee{ \frac{\tau_\la( y)}{h_D(y)}}} &\sleq h_D^{-2}(y)\eee{\tau_\la(y) h_D'(y) +\tau_\la'(y) h_D(y)}\\
&\sleq h_D^{-2} \la (\la y)^{(1-d)/2} e^{\la y} \eee{(y-1)h_D'(y)+ h_D(y)} \\
&\sleq  h_D^{-1}(y) \la (\la y)^{(1-d)/2} e^{\la y} .
}
Therefore,
\spx{
J_3 
&\sleq  h_D^{-1}(y) \int_{1/y}^\8 \la^{d} (\la x)^{(1-d)/2} e^{-\la x} (\la y)^{(1-d)/2} e^{\la y}\, d\la \\
&\sleq h_D^{-1}(y) (xy)^{(1-d)/2} \int_{1/y}^\8  \la e^{-\la(x-y)} \, d\la \\
&\sleq h_D^{-1}(y) (xy)^{(1-d)/2} (x-y)^{-2} \int_{(x-y)/y}^\8  \la e^{-\la} \, d\la\\
&\sleq h_D^{-1}(y) (xy)^{(1-d)/2} e^{-c(x-y)/y} (x-y)^{-2}\\
&\sleq h_D^{-1}(y) \frac{x^{1-d}}{|x-y|^2},
}
where in the last inequality we have used the fact that $e^{-c(x-y)/y} \sleq (x/y)^{(1-d)/2}.$ This ends the proof of \eqref{Riesz_est_D} and Proposition \ref{prop_kernel}.
\end{proof}

\subsection{Proof of Theorem \ref{thm_riesz_halfline}}\label{ssec24}

\begin{proof}
We shall only  prove part {\bf 2.} concerning the operator with the Dirichlet boundary condition. The proof will be based on the estimate \eqref{Riesz_est_D}. The proof of part {\bf 1.} is similar, simpler, and based on \eqref{Riesz_est_N}. The details are left to the reader.

By a standard continuity argument it is enough to show that
$$\norm{R_{D} a}_{L^1(X_+,\mu)} \leq C,$$
for any $(X_+, h_D)$-atom $a$. Suppose then that $B = B(y_0, r)$ is a ball (interval) in $X_+ = (0,\8)$ and a function $a$ and $B$ satisfy \eqref{atom_Dir}.
Denote $5B = B(y_0, 5r)$. Using the boundedness of $R_{D}$ on $L^2(X_+, \mu)$ and the doubling condition for $\mu$ we have
$$\norm{R_{D} a}_{L^1(5B, \mu)} \sleq  \mu(5B)^{1/2} \norm{R_{D} a}_{L^2(5B, \mu)} \sleq \mu(B)^{1/2} \norm{ a}_{L^2(5B, \mu)} \leq C.
$$
Let $z = y_0 +r$ be the right end of the interval $B$. 
For $x \in X_+ \setminus 5B$ using $\int a h_D \, d\mu = 0$ we obtain
\spx{
\abs{R_{D} a(x)} &= \abs{\int_B \eee{ \frac{R_{D}(x,y)}{h_D(y)} - \frac{R_{D}(x,z)}{h_D(z)}  } a(y) h_D(y) \,d\mu(y)}\\
&\leq \norm{a}_{L^1(\mu)} \sup_{y\in B}\underbrace{ \eee{h_D(y) \abs{ \frac{R_{D}(x,y)}{h_D(y)} - \frac{R_{D}(x,z)}{h_D(z)}  }}}_{ W(x,y)}.
}
Using \eqref{Riesz_est_D} and the mean value theorem,
\spx{
W(x,y) \sleq h_D(y) |y-z| h_D(\xi)^{-1} |x-\xi|^{-2} x^{1-d}
}
for some $\xi \in (y,z)$. Notice that $h_D(y) \leq h_D(\xi)$. For $y\in B$ and $x\not \in 5B$ we have $|y-z|\leq 2r$ and $|x-\xi| \simeq |x-y_0|$. Thus,
$$R_Da(x) \sleq \sup_{y\in B} W(x,y) \sleq r  |x-y_0|^{-2} x^{1-d}$$
and
$$\norm{R_Da}_{L^1(\mu, (5B)^c} \sleq \int_{(5B)^c} r x^{1-d} |x-y_0|^{-2}\, d\mu(x) \sleq r \int_{|x-y_0|>5r} |x-y_0|^{-2}\, dx\leq C.$$
The proof is finished.

\end{proof}


\section{The Bessel operator on $\wt{X}$.}\label{sec3}
In this section we study the Bessel operator \eqref{Bessel} on the space $\wt{X}=(-\8,-1]\cup[1,\8)$. We prove here  Proposition \ref{thmA} and  Theorems \ref{thmB}, and \ref{thmC}. Recall that $d>2$ and $d\mu(x) = |x|^{d-1}\, dx$. On $\wt{X}$ we consider the measure $\mu$ and the metric $\rho$ given by the formula
\spx{
\rho(x,y) = 
\begin{cases}
|x-y| &  \text{ if }xy>0\\
|x-y|-2 & \text{ if } xy<0
\end{cases}, \qquad x,y\in \wt{X}.\\
}
In other words, the points $-1$ and $1$ are glued together and the metric is Euclidean on each half-line. As a metric space this is obviously isometric to $\RR$, but it is more convenient for us to use this notation. Recall that
\eq{
\label{measure}
\mu(B(x,r)) \simeq r(|x|+r)^{d-1}.
}

The Bessel operator $\wt{L}$ on $\wt{X}$ is defined  by a quadratic form with the core $C^1(\wt{X})$, similarly as in Definition \ref{def_1}, see also \cite[Section 2.4]{Sikora_Hassell_1D}.

For a function $f$  on $\wt{X}$ let $f_e(x) = (f(x) + f(-x))/2$ and $f_o(x) = f(x) - f_e(x)$ be the even and odd components of $f$. Similarly, for a function $g$ on $X_+$ let $g_{even}$ (resp. $g_{odd}$) be the even (resp. odd) extension of $g$ to $\wt{X}$. Notice that for $f,g \in C^1_c(\wt{X})$ we have
$$Q(f,g) = Q(f_e,g_e) + Q(f_o,g_o).$$
Moreover, $f_e|_{X_+}, g_e|_{X_+} \in \Sigma_N$, since $f_e'(1) = g_e'(1) = 0$, see Definition \ref{def_1}. Thus
$$Q(f_e,g_e) = 2 Q_N (f_e|_{X_+}, g_e|_{X_+})$$
and the operator $\wt{L}$ on even functions corresponds to the operator $L_N$ on $X_+$ in the sense that
$\wt{L} (f_e) = \eee{L_N(f_e|_{X_+})}_{even}.$ Similarly, for an odd function $f\in C_c^1(\wt{X})$ we have $f(1)=0$ and
$\wt{L} (f_o) = \eee{L_D(f_o|_{X_+})}_{odd}.$
Therefore
\begin{equation} \label{SymmetricDecomposition}
\wt{L}(f) = \eee{L_N(f_e|_{X_+})}_{even} + \eee{L_D(f_o|_{X_+})}_{odd}.
\end{equation}

The Bessel operator on $\wt{X}$ was studied in \cite{Sikora_Hassell_1D} and we refer the reader to this paper for the precise definition of $\wt{L}$. As we mentioned in the introduction, it is a model that is related to studying Laplace-Beltrami operators on the manifolds with ends of the form $\RR^d \# \RR^d$.

The analysis of $\wt{L}$ is particularly interesting because of a phenomenon that this operator has two linearly independent bounded harmonic functions. These harmonic functions can be defined in terms of probability in the following way. Let $h_+$ (resp. $h_-$) be the probability that a process starting at the point $x\in \wt{X}$ exits to infinity via ''the right (resp. left) end''. Such functions are harmonic and bounded. In fact we have
$$
h_+(x) = 
\begin{cases}
1-x^{2-d}/2, & x\in X_+\\
|x|^{2-d}/2, & x\in X_-.
\end{cases}
, \quad h_-(x) =
\begin{cases}
x^{2-d}/2, & x\in X_+\\
1-|x|^{2-d}/2, & x\in X_-
\end{cases}
.
$$
Denote $\wh{h}_N \equiv 1$ on $\wt{X}$ and let $\wh{h}_D$ be the odd extension of $h_D$ given by \eqref{eq_hD}. Notice that
$$h_++h_- = \wh{h}_N, \qquad h_+ - h_- = \wh{h}_D.$$
The next lemma gives a crucial relation between kernels of the semigroups related to $\wt{L}$, $L_N$ and $L_D$. Using this relation the proofs of Proposition \ref{thmA} and Theorems \ref{thmB} and \ref{thmC} boil down to the analysis of $L_N$ and $L_D$ on $X_+$.

\lem{lem2}{
Let $d>2$. The kernel $T_t(x,y)$ of $\exp(-t\wt{L})$ can be expressed in the following way:
\alx{
T_t (x,y) &=\frac{1}{2} \eee{ T_{t,N} (x,y) + T_{t,D} (x,y)}, && x,y \in X_+\\
T_t(x,y) &= \frac{1}{2} \eee{T_{t,N} (x,-y) - T_{t,D} (x,-y)}, &&x \in X_+, y\in X_-.
}
For $x\in X_-$ we have $T_t(x,y) = T_t(-x,-y)$,}

\pr{
From \eqref{SymmetricDecomposition} we deduce that the semigroup $T_t = \exp(-t\wt{L})$ is given by
\spx{
T_t f(x) = \eee{T_{t,N} (f_e|_{X_+})}_{even}(x) +
\eee{T_{t,D} (f_o|_{X_+})}_{odd}(x)
}
and Lemma \ref{lem2} follows from this relation.
}

\subsection{Proof of Proposition \ref{thmA}} \label{ssec3.3}
\pr{
To prove Proposition \ref{thmA} recall that on $X_+$ the operator $L_N$ satisfies \eqref{gauss1}, whereas the operator $L_D$ satisfies \eqref{gauss2} with the function $h_D(x) = 1-x^{2-d}$. We claim that $h_D$ belongs to the Reverse H\"older class $RH_\8(X_+)$, i.e. for every ball (interval) $B$ in $X_+$ we have:
\eq{\label{RH}
\sup_{y\in B} h_D(y) \leq C \mu(B)^{-1} \int_B h_D(x)\, d\mu(x),
}
where $C$ does not depend on $B$. Indeed, let $[a,b] \subseteq [1,\8)$ and denote $c=(a+b)/2$. By the doubling condition of $\mu$ we have
$$\mu(B)^{-1} \int_a^b h_D(x) \, d\mu(x) \sgeq \mu([c,b])^{-1} \int_c^b h_D(x) \, d\mu(x) \geq h_D(c).$$
Since $h_D(x) \simeq (x-1)/x$, we have 
$$h_D(c) \geq h_D((b+1)/2) \simeq \frac{b-1}{b+1} \geq  \frac{b-1}{b} \simeq h_D(b) = \sup_{y\in B} h_D(y)$$
and the claim is proved.

From \cite[Theorem 4.3]{Preisner_Sikora} we have  the space $H^1_{\wt{L}}(\wt{X})$ is equal to $H^1_{at}(\mu, \wh{h}_N, \wh{h}_D)$, where $\wh{h}_N \equiv 1$ on $\wt{X}$ and $\wh{h}_D$ is the odd extension of $h_D$ to $\wt{X}$. This gives Proposition \ref{thmA} since
$$h_+(x) = \frac{1}{2}\eee{\wh{h}_N(x) +\wh{h}_D(x)} \quad \text{and} \quad  h_-(x) = \frac{1}{2}\eee{\wh{h}_N(x) -\wh{h}_D(x)}$$
and, therefore, the spaces $H^1_{at}(\mu, \wh{h}_N, \wh{h}_D)$ and $H^1_{at}(\mu, h_+, h_-)$ are equal.
}


\subsection{Proof of Theorem \ref{thmB}}\label{ssec3.4}
\pr{
Recall that the operator $\wt{L}$ can be decomposed in terms of $L_N$ and $L_D$, see \eqref{SymmetricDecomposition}. Therefore,
\sp{\label{eq_Riesz_decomp}
\wt{R}f = (R_N(f_e|_{X_+}))_{even} + 
(R_D(f_o|_{X_+}))_{odd}.
} 

From Proposition \ref{thmA} we know that $H^1_{\wt{L}}(\wt{X})$ coincide with $H^1_{at}(\mu, h_+, h_-)=H^1_{at}(\mu, \wh{h}_N, \wh{h}_D)$, where $\wh{h}_N \equiv 1$ on $\wh{X}$ and $\wh{h}_D$ is the odd extension of $h_D$ to $\wh{X}$, see the proof of Proposition  \ref{thmA}. Similarly as in \ref{ssec24}, by a standard continuity argument, it is enough to show that
\[
\|\wt{R}a\|_{L^1(\wt{X},\mu)} \leq C,
\]
for any $(\mu,\wh{h}_N,\wh{h}_D)$-atom $a$ with a constant $C>0$ that does not depend on $a$. Assume then that $a$ is such atom and denote by $a_e$ and $a_o$ the even and odd part of $a$, respectively. Denote, $b_e := a_e|_{X_+}$, $b_o := a_o|_{X_+}$. Notice that
\eq{\label{eq_norm_decomp}
\|a\|_{L^2(\wt{X},\mu)} \simeq 
\|a_e\|_{L^2(\wt{X},\mu)} + \|a_o\|_{L^2(\wt{X},\mu)} = 2\eee{\|b_e\|_{L^2(X_+,\mu)} + \|b_o\|_{L^2(X_+,\mu)}}.
}
Moreover,
$$0 = \int_{\wt{X}} a(x) \, d\mu(x) = \int_{\wt{X}} a_e(x)\, d\mu(x) = 2 \int_{X_+} b_e(x)\, d\mu(x),$$
and, therefore, $b_e$ is an $(\mu,h_N)$-atom on $(X_+, \mu)$. Similarly,
$$0 = \int_{\wt{X}} a(x) \wh{h}_D(x) \, d\mu(x) = \int_{\wt{X}} a_o(x) \wh{h}_D(x)\, d\mu(x) = 2 \int_{X_+} b_o(x) h_D(x)\, d\mu(x),$$
and $b_o$ is an $(\mu,h_D)$-atom on $(X_+, \mu)$.

Using \eqref{eq_Riesz_decomp} and Theorem \ref{thm_riesz_halfline} we get
\spx{
\|\wt{R}a\|_{L^1(\wt{X},\mu)} 
&\leq
2||R_N b_e||_{L^1(X_+,\mu)} +
2||R_D b_o||_{L^1(X_+,\mu)} \leq C.
}
}

\begin{Rem}\label{remark}
    Let $d>2$ and set $g(x) = \chi_{[2,3]} (x) - \chi_{[-3,-2]}(x)$. Then $g$ is in the Coifman-Weiss atomic Hardy space $H^1_{at}(\wt{X})$ on $\wt{X}$, c.f. \eqref{Coifman-Weiss}, but $\wt{R}g \not\in L^1(\wt{X}, \mu)$.
\end{Rem}
\begin{proof}
    It is clear that $g$, after multiplying by a certain constant, is a Coifman-Weiss atom on $\wt{X}$, c.f. \eqref{Coifman-Weiss}. Notice that $\wt{R} g$ is an odd function and for $x\in X_+$ it is equal to $R_D (\chi_{[2,3]})$. By Proposition \ref{prop_Riesz_kernel} and Lemma \ref{tau_lemma} for $y\in[2,3]$ and $x\geq 4$ we get $R_D(x,y)$ is non-positive and
    \spx{
    -R_D(x,y) \geq \int_0^{1/x} \la^{d-1} (\la x)^{1-d}\, d\la \sgeq x^{-d}.
    }
Therefore, $\wt{R}g \not\in L^1(\wt{X},\mu)$.
\end{proof}

\subsection{Proof of Theorem \ref{thmC}}

\pr{ 
Consider first $x$ and $y$ on the same side. Without loss of generality let $x,y >0$. By Lemma \ref{lem2} we have that $T_t(x,y) = \frac{1}{2}\eee{T_{t,N}(x,y) +T_{t,D}(x,y)}$. Then \eqref{ULmixedBounds0} follows from \eqref{gauss1} and the fact that the Dirichlet heat kernel is dominated by the Neumann heat kernel.

Now,  consider $x$ and $y$ such that $xy<0$. We have  $T_t(x,y) = \frac{1}{2}\eee{T_{t,N}(|x|,|y|) -T_{t,D}(|x|,|y|)}$. We shall prove \eqref{ULmixedBounds} in a few steps.

{\bf Step 1.} Assume first that at least one of the points $x,y$ is in the middle past, i.e. in $[-2,-1]\cup [1,2]$. Without loss of generality we may assume that $y\in [-2,-1]$ and $x\geq 1$. Let $\widehat{T}_t$ be the process associate with the Bessel operator $\widehat{L}$ on $\widehat{X} :=\wt{X} \cap [-3,\8)$ with the Dirichlet boundary condition at $x=-3$, c.f. Definition \ref{def_1}. Then,
\eq{\label{dwa_procesy}
T_t(x,y) \geq \widehat{T}_t(x,y), \quad x,y \in \widehat{X},
}
since the process $\widehat{T}_t$ can be considered as the process $T_t$ up to the first arrival to the boundary point $x=-3$. Let $\widehat{h}$ be the harmonic function associated with $\widehat{T}_t$. Similarly as in Section \ref{ssec22}, using \cite{Gyrya_Saloff-Coste}, we have that
\eq{\label{gauss3}
\mu_{\widehat{h}^2}(B(x,\st))^{-1} \exp\eee{-\frac{|x-y|^2}{c_1 t}}\sleq \frac{\widehat{T}_{t}(x,y)}{\widehat{h}(x)\widehat{h}(y)} \sleq  \mu_{\widehat{h}^2}(B(x,\st))^{-1} \exp\eee{-\frac{|x-y|^2}{c_2 t}}.
}
The harmonic function $\widehat{h}(x)$ can be computed by finding the $C^1$ solution of the equation $Bf = 0$ on $\widehat{X}$ with the boundary condition $\widehat{h}(-3)=0$. It follows that
\eq{
\label{comp1}
\widehat{h}(x) \simeq 1, \qquad x\in[-2,-1]\cup [1,\8).
}
Then, \eqref{ULmixedBounds} for $x\geq 1$ and $y\in[-2,-1]$ follows from \eqref{dwa_procesy}, \eqref{gauss3} and \eqref{comp1}.

{\bf Step 2a.}
For $x,y\geq 1$ and $t,c_1,c_2>0$ let us denote 
\[
\Pi_{c_1,c_2}(x,y,t) =
 \int^t_0 \frac{(t-s)^{-1/2}}{(x + \sqrt{t-s})^{d-1}} \frac{y-1}{y}
\frac{s^{-3/2}}{(s + y)^{\frac{d-3}{2}}}
\exp\eee{-\frac{(x-1)^2}{c_1(t-s)}-\frac{(y-1)^2}{c_2s}}  ds.
\]
We claim that the kernel $T_t(x,y)$ satisfies
\begin{equation}
\label{mixedKernelBound}
\Pi_{c_1',c_2'}(x,y,t) \sleq 
T_t(x,-y) \sleq
\Pi_{c_1'',c_2''}(x,y,t),
\end{equation}
for some $c_1',c_2',c_1'',c_2''>0$ and $x,y \geq 1$.

Consider $x,-y\in \wt{X}$ such that $x,y\geq 1$. By Lemma \ref{lem2} we have that
$$T_t(x,-y) = \frac{1}{2}\eee{T_{t,N}(x,y) - T_{t,D}(x,y)}.$$
To estimate the difference above, we shall use a  probabilistic argument. Recall that with the operators $L_N$ and $L_D$ (or with the related quadratic forms), one can associate the unique Hunt processes $X_t^N$ and $X_t^D$ on $X_+ = [1,\8)$, where the latter one is the process that is killed at the boundary $x=1$. The paths of these processes are continuous by \cite[Theorem 7.2.2]{DirichletForms}. Denote by $\tau_1$ the stopping time
$\tau_1 =\inf\{t>0: X_t^N = 1\}=\inf\{t>0: X_t^D = 1\}.$ A well-known Dynkin-Hunt formula, see e.g. \cite[p.42]{Gyrya_Saloff-Coste}, states that
\[
T_{t,N}(x,y) -T_{t,D}(x,y) = \EE^{-y}[{\bf 1}_{\{\tau_1\leq t\}}T_{t-\tau_1,N}(x,1)]=\int_0^t T_{t-s,N}(x,1) \, d\PP^{-y}(\tau_1 \leq s),
\]
compare \cite[Theorem 2.3]{Grigoryan_Ishiwata_Saloff-Coste}.
This observation allows us to use \cite[Theorem 2]{HittingTimes} to estimate the integral above by
\begin{equation}
\label{eq435}    
T_t(x,-y) \simeq \int_0^t T_{t-s,N}(x,1)\frac{(y-1)y^{d-3}}{1+y^{d-2}}\frac{\exp\eee{-\frac{(y-1)^2}{4s}}}{(s + y)^{\frac{d-3}{2}}}  \frac{ds}{s^{3/2}}.
\end{equation}

The estimate \eqref{mixedKernelBound} follows directly from \eqref{eq435} and \eqref{gauss1}.

{\bf Step 2b.} What is left to prove is \eqref{ULmixedBounds} in the case $|x|,|y|\geq 2$. Without loss of generality we can consider the kernel $T_t(x,-y)$ for $x\geq y\geq 2$. In this case $\rho(x,y) \simeq x$. Let $c_1, c_2>0$ be fixed constants. Using \eqref{mixedKernelBound}, it is enough to prove that for $x\geq y\geq 2$ and $t>0$ we have
\eq{
\label{est321}
\frac{ y^{2-d}}{\mu(B(x,\sqrt{t}))} \exp\eee{-\frac{x^2}{c_1' t}} \sleq \Pi_{c_1,c_2}(x,y,t) \sleq \frac{ y^{d-2}}{\mu(B(x,\sqrt{t}))} \exp\eee{-\frac{x^2}{c_1'' t}},
}
with some $c_1',c_1''>0$. Notice that $x \simeq x-1$ and $y\simeq y-1$. Write
\spx{
\int^t_0 \frac{(t-s)^{-1/2}}{(x + \sqrt{t-s})^{d-1}}\frac{s^{-3/2}}{s^{\frac{d-3}{2}} + y^{\frac{d-3}{2}}}
 \exp\eee{-\frac{x^2}{b_1(t-s)}-\frac{y^2}{b_2s}} ds
= \underbrace{\int_0^{\frac{t}{2}}...\, ds}_{\Pi^{[1]}_{b_1,b_2}(x,y,t)}+\underbrace{\int_{\frac{t}{2}}^t...\, ds}_{\Pi^{[2]}_{b_1,b_2}(x,y,t)}.
}
Then $\Pi_{c_1,c_2}(x,y,t) \sleq \Pi^{[1]}_{4c_1,4c_2}(x,y,t) + \Pi^{[2]}_{4c_1,4c_2}(x,y,t)$ and $\Pi_{c_1,c_2}(x,y,t) \sgeq \Pi^{[1]}_{c_1,c_2}(x,y,t)$.
Thus, it is enough to prove that
\al{
\label{ff1}
\Pi^{[1]}_{4c_1,4c_2}(x,y,t) &\sleq  \frac{ y^{d-2}}{\mu(B(x,\sqrt{t}))} \exp\eee{-\frac{x^2}{c_1'' t}},\\
\label{ff2}
\Pi^{[1]}_{c_1,c_2}(x,y,t) &\sgeq \frac{ y^{2-d}}{\mu(B(x,\sqrt{t}))} \exp\eee{-\frac{x^2}{c_1' t}},\\
\label{ff3}
\Pi^{[2]}_{4c_1,4c_2}(x,y,t) &\sleq  \frac{ y^{d-2}}{\mu(B(x,\sqrt{t}))} \exp\eee{-\frac{x^2}{c_1'' t}}.
}

Let us start by considering $\Pi^{[1]}$. We have that $t-s\simeq t$ for $s\in[0,t/2]$. Using \eqref{measure},
\spx{\label{Pi1}
\Pi^{[1]}_{4c_1,4c_2}(x,y,t) 
&\sleq \frac{1}{\mu(B(x,\sqrt{t}))}\exp\eee{-\frac{x^2}{c_3 t}}\underbrace{\int_0^{\frac{t}{2}} \frac{s^{-3/2}}{(s+y)^{\frac{d-3}{2}}} \exp\eee{-\frac{y^2}{c_4 s}} ds}_{\Psi(y,t)}, \\
}
with some constants $c_3,c_4>0$. Similarly,
$$\Pi^{[1]}_{c_1,c_2}(x,y,t) \sgeq \frac{1}{\mu(B(x,\sqrt{t}))}\exp\eee{-\frac{x^2}{c_3' t}}\Psi(y,t),
$$
with a possibly different constant $c_4$ in $\Psi$.

We claim that
\eq{
\label{psi_est}
y^{2-d} \exp\eee{-\frac{y^2}{c_5 t}} \sleq \Psi(y,t) \sleq y^{2-d}
}
with some $c_5=c_4/4$. Observe that \eqref{ff1} and \eqref{ff2} follow from \eqref{psi_est}, since $x\geq y$.

Now, we prove the claim \eqref{psi_est}. Let us consider three cases.

{\bf Case 1:}
For $\frac{t}{4} \leq y$
\alx{
\Psi(y,t)
&\sgeq \int_{\frac{t}{4}}^{\frac{t}{2}} \frac{s^{-3/2}}{(s+y)^{\frac{d-3}{2}}}\exp\eee{-\frac{y^2}{c_4s}} ds
\sgeq \frac{t^{-1/2}}{y^{\frac{d-3}{2}}}\exp\eee{-\frac{y^2}{c_5t}} \\
&\sgeq y^{\frac{2-d}{2}}\exp\eee{-\frac{4y^2}{c_4t}}
\sgeq y^{2-d}\exp\eee{-\frac{y^2}{c_5t}},\\
\Psi(y,t)
&\sleq \int_0^y s^{-3/2} y^{\frac{3-d}{2}} \eee{\frac{s}{y^2}}^{\frac{d}{2}+1}\, ds \simeq y^{-d}\sleq y^{2-d}.
}

{\bf Case 2:}
Next, if $y < \frac{t}{4} < y^2$, then
\spx{
\Psi(y,t) 
&\sgeq \int_{\frac{t}{4}}^{\frac{t}{2}} s^{-d/2}\exp\eee{-\frac{y^2}{c_4s}}\, ds
\sgeq t^{1-d/2}\exp\eee{-\frac{4y^2}{c_4t}}
\sgeq y^{2-d}\exp\eee{-\frac{y^2}{c_5t}}, \\
\Psi(y,t) 
&\sleq \int^y_0 ... + \int^{\frac{t}{2}}_y ... 
\sleq \int^y_0 \frac{s^{-3/2}}{y^{\frac{d-3}{2}}} \eee{\frac{s}{y^2}}^{d/2+1} ds 
+ \int_0^{\infty} s^{-d/2}\exp\eee{-\frac{y^2}{c_4s}} \, ds \\
&\simeq y^{-d} + y^{2-d}
\leq y^{2-d}.
}

{\bf Case 3:} $y^2 < t/4$. Now
\spx{
\Psi(y,t)
&\geq \int^{2y^2}_{y^2} ...
\simeq \int^{2y^2}_{y^2} s^{-d/2} ds
\simeq y^{2-d}
\simeq y^{2-d}\exp\eee{-\frac{y^2}{c_5t}}, \\
\Psi(y,t) 
&\sleq \int^y_0 ... + \int^{\frac{t}{2}}_y ... 
\sleq \int^y_0 \frac{s^{-3/2}}{y^{\frac{d-3}{2}}} \eee{\frac{s}{y^2}}^{d/2+1} ds 
+ \int_0^{\infty} s^{-d/2}\exp\eee{-\frac{y^2}{c_4s}} \,ds \\
&\simeq y^{-d} + y^{2-d}
\leq y^{2-d},
}
similarly as in {\bf Case 2.} The claim \eqref{psi_est} is proved.


We now consider $\Pi_{4c_1,4c_2}^{[2]}$. Write
\spx{
\Pi_{4c_1,4c_2}^{[2]}
(x,y,t)
&\sleq \exp\eee{-\frac{y^2}{4c_2 t}}\frac{t^{-\frac{3}{2}}}{(t+y)^{\frac{d-3}{2}}}\int_{\frac{t}{2}}^t \frac{(t-s)^{-1/2}}{(x+\sqrt{t-s})^{d-1}}e^{-\frac{x^2}{c_1(t-s)}}ds\\
&= \exp\eee{-\frac{y^2}{4c_2 t}}\frac{t^{-\frac{3}{2}}}{(t+y)^{\frac{d-3}{2}}}\underbrace{\int_0^{\frac{t}{2}} \frac{s^{-1/2}}{(x+\sqrt{s})^{d-1}}e^{-\frac{x^2}{c_1s}}ds}_{\Phi(x,t)}.
}
For $\frac{t}{4}<x^2$ we have
\spx{
\Phi(x,y,t)
\simeq \int_0^{\frac{t}{2}}s^{-\frac{1}{2}}x^{1-d}e^{-\frac{x^2}{c_1s}}ds \sleq  x^{1-d} \exp\eee{-x^2/c_1 t} \int_0^{x^2} s^{-1/2}\, ds\simeq x^{2-d} \exp\eee{-\frac{x^2}{c_1 t}}.
}
If $\frac{t}{4} \ge x^2$, then
\spx{
\Phi(x,y,t) \sleq \int_0^\8 s^{-d/2} \exp\eee{-\frac{x^2}{c_1s}}\, ds =x^{2-d} \int_0^\8 s^{-d/2} \exp\eee{-\frac{1}{c_1s}}\, ds \simeq x^{2-d} \exp\eee{-\frac{x^2}{c_1t}}.
}
This gives
\spx{
\Pi^{[2]}_{4c_1,4c_2}(x,y,t)
&\sleq \exp\eee{-\frac{y^2}{4c_2 t}}\frac{t^{-\frac{3}{2}}}{(t+y)^{\frac{d-3}{2}}} x^{2-d}e^{-\frac{x^2}{c_1t}} \\
&\sleq \underbrace{\exp\eee{-\frac{y^2}{4c_2 t}}\frac{t^{-1/2}}{(t+y)^{\frac{d-3}{2}}}y^{d-2}}_{G(t,y)} y^{2-d}  t^{-1/2} x^{1-d} \frac{x}{\sqrt{t}}e^{-\frac{x^2}{c_1t}}\\
&\sleq y^{2-d} \mu(B(x,\sqrt{t}))^{-1} e^{-\frac{x^2}{2c_1t}},
}
where we have used \eqref{measure} and the fact that $G$ is a bounded function.
}

{\bf Funding.}
M.P and A.S. were partly supported by Australian Research Council (ARC) Discovery Grant   DP200101065.


\def\cprime{$'$}
\providecommand{\bysame}{\leavevmode\hbox to3em{\hrulefill}\thinspace}
\providecommand{\MR}{\relax\ifhmode\unskip\space\fi MR }
\providecommand{\MRhref}[2]{%
  \href{http://www.ams.org/mathscinet-getitem?mr=#1}{#2}
}
\providecommand{\href}[2]{#2}

\end{document}